\documentclass{ifacconf}

\newtheorem{theorem}{Theorem}[section]

\newtheorem{definition}[theorem]{Definition}
\newtheorem{remark}[theorem]{Remark}
\newcommand{\R}{\mathbb{R}}

\newcommand{\B}{\mathbb{B}}
\newcommand{\hsm}{\hspace{-0.4mm}}
\newcommand{\hsmm}{\hspace{-0.2mm}}
\newcommand{\vsm}{\vspace{-0.5mm}}
\usepackage{amsfonts}
\usepackage{comment}
\usepackage{amsmath}
\usepackage{xcolor}
\usepackage{graphicx}      
\usepackage{natbib}        
\begin{document}
\begin{frontmatter}

\title{A Sweeping Process Control Problem Subject To Mixed Constraints} 

\thanks[footnoteinfo]{Nathalie T. Khalil  acknowledges the support of ARISE AL under Grant
LA/P/0112/2020; of SYSTEC under Grant UIDB/00147/2020 and Grant UIDP/00147/2020; and of projects SNAP under Grant NORTE-01-0145-FEDER-000085 and MAGIC under Grant PTDC/EEIAUT/
32485/2017, funded by COMPETE2020-POCI and FCT/MCTES. This work had also the partial support of the Department of Mathematics of the Universidad Autónoma Metropolitana-Iztapalapa, Mexico.}

\author[First]{Karla L. Cortez} 
\author[Second]{Nathalie T. Khalil} 
\author[First]{Julio E. Solís}

\address[First]{Department of Mathematics, Universidad Autónoma Metropolitana-Iztapalapa, Mexico (e-mail:  kcortez@xanum.uam.mx; jesd@xanum.uam.mx).}
\address[Second]{Department of Electrical and Computer Engineering, University of Porto, Portugal (e-mail: nathalie@fe.up.pt)}

\begin{abstract}                
In this study, we investigate optimal control problems that involve sweeping processes with a drift term and mixed inequality constraints. Our goal is to establish necessary optimality conditions for these problems. We address the challenges that arise due to the combination of sweeping processes and inequality mixed constraints in two contexts: regular and non-regular. This requires working with different types of multipliers, such as finite positive Radon measures for the sweeping term and integrable functions for regular mixed constraints. For non-regular mixed constraints, the multipliers correspond to purely finitely additive set functions.


\end{abstract}

\end{frontmatter}

\section{Introduction}\label{sec:introduction}
The goal of this study is to examine an optimal control problem that features a sweeping process with a drift term, subject to a mixed state/control constraint condition. Specifically, we are focused on the problem over a fixed time interval of $[0,1]$ which is formulated as follows:
	\begin{eqnarray}
		(P)\mbox{ Minimize} & & g(x(1))\nonumber\\ \nonumber
		\mbox{s.t.} & &   \dot x(t) \in f(x(t),u(t)) - N_C(x(t)) \;\;[0,1]\textrm{-a.e.}\nonumber \\
		&&  x(0) \in C_0 \subset C\;\; \nonumber \\
		&&  h(x(t),u(t))\leq 0,\;  
  [0,1]\textrm{-a.e.}\label{mixed constraint}
	\end{eqnarray}
	Here, $g:\R^n \to \R$, $f: \R^n\times \R^m \to \R^n$, and $h:\R^n \times \R^m \to \R$ are the objective function, the drift term of the sweeping process, and the inequality mixed (state/control) constraint, 
	respectively. $N_{A}(z)$ is the limiting normal cone to the closed set $A$ at point $z$, in the sense of \cite{mordukhovich2006variational}, representing the sweeping process. $C_0 \subset \R^{n}$ the left-endpoint constraint set, and $C$ is the sweeping set specified by
	\begin{equation}\label{definition sweeping set}
	    C:=\{ x \in \R^n:\psi(x) \le 0 \}
	\end{equation} 
	where $\psi: \R^n \to \R$, such that $0\in \text{int }C.$ (Here \text{int }$C$ represents the interior of the set $C$.)	For the sake of simplicity, we consider $\psi$ to be convex and smooth, and in this context, $N_{C}(x)$ is the normal cone to the convex set $C$ at $x$, which will take the simplified form: 
\begin{equation}  \label{def: normal cone}
N_C(x):= \begin{cases}
	\{ v\in\R^n|\langle v,\bar x- x\rangle \leq 0\,\forall \,\bar x\in C\} & \textrm{if }  x\in C 
	
		\\  \emptyset & \textrm{otherwise}.
\end{cases}
\end{equation}
 Moreover, in this paper, we refrain from imposing any constraint on the control set, denoted $U$. In fact, we consider it as the whole space $U:=\R^m$, and the controls $u(.)$ to be essentially bounded with values in $\R^m$, i.e. $u\in L^{\infty}([0,1],\R^m)$.

The purpose of this article is to establish the necessary optimality conditions for (P) in two different scenarios: regular and non-regular mixed constraint. Due to the challenging nature of the sweeping process term and the mixed constraints, this formulation has not been previously considered. The obstacles to overcome include addressing a non-Lipschitz differential inclusion with a discontinuous normal cone and managing several multipliers associated with the sweeping term and the mixed constraint type (regular or non-regular).


Mixed-constrained optimal control problems for ordinary differential equations have been extensively studied in the literature, with the necessary optimality conditions in Dubovitskii-Milyutin form explored in works such as \cite{makowski1974optimal}, \cite{dmitruk2009development}, and \cite{clarke2010optimal} for the regular case. The optimality conditions in Gamkrelidze's form have been established in various works, including \cite{arutyunov2010maximum}, \cite{arutyunov2016investigation}, and \cite{khalilpereiraACC2021}, where the latter studied a differential inclusion dynamic with an extra pure state constraint. More recently, the non-regular case has been investigated in \cite{Becerriletal2021}, \cite{becerriletal2022}, and \cite{dmitruk2022local}.

Regular and non-regular mixed-constrained problems differ in two aspects: the quality of multipliers and the information provided by the necessary optimality conditions. Regular problems have an integrable function as the multiplier associated with the mixed constraint, and the Maximum Principle gives complete information. In contrast, non-regular problems require considering the dual space of essentially bounded functions, $(L^\infty)^*$, which includes ``charges'' as new multipliers associated with the mixed constraints. Additionally, the Maximum Principle does not provide the same level of information as in regular problems due to the absence of the maximization condition.


    
Control problems involving sweeping processes have been also widely studied in the literature. The concept of sweeping processes was introduced by Jean-Jacques Moreau for problems on elastoplasticity, cf. \cite{moreau1976application}. The main focus of research includes the well-posedness, existence, regularity of solutions, and numerical aspects of sweeping, cf. \cite{moreau1999numerical}, \cite{maury2008mathematical}, \cite{venel2011numerical}, \cite{brokate2013optimal}. Various necessary optimality conditions have been established for different types of sweeping sets, including moving hyperplanes, moving convex polyhedra, and uniform prox-regular nonconvex moving sets, as in \cite{colombo2012optimal}, \cite{colombo2015discrete}, and \cite{cao2019optimal} respectively. Recently, necessary conditions were derived for a constant, convex sweeping set with the addition of a pure state constraint condition in \cite{khalil2022maximum}.
	
This paper studies problem $(P)$, which involves a constant sweeping set $C$, smooth data, and a scalar inequality mixed constraint. The mixed constraint is initially assumed to be regular (see Definition \ref{def:regular-mixed traj}), and necessary conditions for problem $(P)$ are established based on a close approach of \cite{khalil2022maximum}.

The problem is then reformulated with the sweeping term expressed as a pure state constraint and a non-regular mixed equality constraint as follows:
\begin{eqnarray}
		&&   \dot x(t) = f(x(t),u(t)) - v(t)\nabla_x \psi (x(t)) \;[0,1]\textrm{-a.e.}\nonumber \\
		&&    \psi(x(t)) \leq 0 \;\; \forall t \in [0,1] \label{state constraint refomulated pb}\\
		&&  v(t)\psi(x(t)) =0 \;[0,1]\textrm{-a.e.}\label{mixed non regular constraint refomulated pb}  \\
		&& (v,u)\in {\cal V} \times \R^m\nonumber
	\end{eqnarray}

	where ${\cal V}=\{v\in L_\infty(\left[0,1\right];\R): v(t)\geq 0,\ \forall t\}$. In this problem, the control variable is the pair $(u, v)$, inequality (\ref{state constraint refomulated pb}) is a pure state constraint condition, while (\ref{mixed non regular constraint refomulated pb}) is a mixed equality constraint, which is {\it non-regular}.
Since most works on maximum principles for mixed constraints and pure state constraints only consider regular mixed constraints, the presence of this non-regular mixed constraint in this reformulation makes it difficult to handle necessary optimality conditions. However, in this paper, we follow the approach of \cite{Becerriletal2021} and \cite{dmitruk2009development} to derive necessary optimality conditions for problem $(P)$, even if
the constraint (\ref{mixed constraint}) is not regular, at the cost of losing the maximum condition and the appearance of purely finitely additive set functions (known as charges) as multipliers.

	
	
	The paper is organized as follows: in the next section, we give the general assumptions on the data of the problem. In section \ref{section:PMP regular mixed constraint}, we establish the necessary optimality conditions of $(P)$ for a regular mixed constraint, while in section \ref{section:PMP non-regular mixed constraint}, the optimality conditions are derived for the case of non-regular mixed constraints. The proof of the  regular case is subject to section \ref{section:proof}. We end the paper with conclusions and future avenues of work.
	
	{\bf Notation.} We denote by $\partial \varphi$ the Mordukhovich (limiting) subdifferential of the function $\varphi$, and by $N_{A}(z)$ the Mordukhovich (limiting) normal cone to the closed set $A$ at point $z$. We refer the reader to \cite{mordukhovich2006variational,vinter2010optimal} for more details on nonsmooth analysis tools. $AC([0,1];\R^n)$ stands for the set of absolutely continuous functions defined on $[0,1]$ with values in $\R^n$, while $BV([0,1];\R^n)$ is the space of functions of bounded variations defined on $[0,1]$, $NBV([0,1];\R^n)$  the
space of functions of bounded variation, right continuous in $(0,1)$ and vanishing at
$1$, and $\|\cdot\|_{TV}$ is the total variation norm. For a smooth function $\varphi(x,u):\R^n\times\R^m\to\R$, we denote by $\nabla_x\varphi$ and $\nabla_u\varphi$ the gradients w.r.t. $x$ and $u$ respectively, while for vector-valued functions, $D_x$, $D_u$ and $D$ denote, respectively, the Jacobian matrix w.r.t. $x$, $u$ and $(x,u).$ 

\section{Assumptions}\label{section: assumptions}
We impose the following assumptions on the data of the problem. Those assumptions can be weakened, however for the sake of simplicity, we shall consider the needed smoothness on the data.
\begin{itemize}
\item[H1] Function $f$ is $\mathcal B \times\mathcal B$-measurable and continuously differentiable in $(x,u)$. Also there exists a constant $M >  0$ such that $|f(x,u)|\leq  M$ and $|D_x f(x,u)|\leq  M$, for all
$( x, u)$.
\item[H2] The set $f(x,\mathbb{R}^m)$ is compact and convex for all $x$.
\item[H3] Function $\psi$ defining the sweeping set $C$ is twice differentiable and convex. Moreover, there exists a constant $\eta >0\mbox{ s.t. } |\nabla_x\psi(x)| > 2\eta \;\; \forall x\in \partial C$, and $\lim\limits_{|x| \to \infty} \psi(x)= \infty.$
\item[H4] Function $h$ is continuously differentiable in $(x,u)$, uniformly bounded, and $\nabla h(\cdot)$ is bounded for all $(x,u)$.
\item[H5] Set $C_0$ is closed and $C_0 \subset C$.
\item[H6] Function $g$ is Lipschitz continuous.
\end{itemize}

	

\section{Regular mixed constraints} \label{section:PMP regular mixed constraint}

We shall first define the set:
\begin{align}\label{def of mixed set}
	\Omega(x):= \{ u \in \R^m \ : \   h(x,u) \le 0 \},
\end{align}
and we give the definition of regularity imposed on the mixed constraints, cf. \cite{arutyunov2016investigation}, which in our context ($h$ scalar) reduces to a simpler definition as shown below. 
\begin{definition}[Regularity]\label{def:regular-mixed traj}   Mixed constraint $h^j$, for $j=1,\ldots,r$ (i.e. vector-valued), is called regular if the linear independence of the set $\{\nabla_u h^j (x,u)\hsm : j \in J(x,u)\}$ holds, where $J(x,u) := \{j\hsm : h^j (x,u)= 0\}$ is the set of active indices.

   When the mixed constraint is scalar (i.e., for $j=1$), the regularity reduces to the following: for any pair $(x(t),u(t))$ satisfying $h(x(t),u(t))=0$, we have $$\nabla_u h(x(t), u(t)) \neq 0.$$
   
\end{definition}

In what follows, $\nabla_x^2\psi$ and $Q(x)$ stand for $D_x(\nabla_x\psi)$ and $ \nabla_x\psi(x)\otimes\nabla_x\psi(x)$, respectively, where $\otimes$ denotes the external product, and for any Borel measure $d\mu$, $S_\mu $ is its support. Define $I_\psi:=\{t\in [0,1]: \psi( x^* (t))=0\}$.

\begin{theorem}\label{theorem PMP regular}
	Let $(x^*,u^*)$ be a solution to $(P)$. Assume that H1-H6 hold and that the mixed constraints are regular. Then, there exists a set of multipliers $(p,\lambda_0,\nu,\eta)$, with $ p\in BV([0,1];\R^n)$, $\lambda_0\ge 0 $, and $\nu: [0,1]\to \R$ measurable and essentially bounded with a non-negative component, a finite positive Radon measure $d\eta$, such that $S_\eta= I_{\psi}$, $\xi \in L^2([0,1],\R)$ verifying $\xi(t)\ge 0$ a.e. $t$, $\xi(t)= 0$ for all $t\in [0,1]\setminus I_{\psi}$, and a positive number $\kappa$  satisfying the following conditions:
\begin{itemize}
	\item[1.] Nontriviality: $ \|p\|+\|\eta\| + \lambda_0\neq 0$.
	\item[2.] Measure-driven adjoint equation:\vsm
	\begin{eqnarray*} &&\hspace{-.7cm} -dp(t) =  p(t) \big[ D_xf(x^*(t), u^*(t))
	-  \xi(t)\nabla_x^2\psi(x^*(t))\big] dt \\
		&& \hspace{.7cm}- p(t)Q(x^*(t)) d\eta(t) -\nu(t)\nabla_x h(x^*(t),u^*(t))dt.
	\end{eqnarray*}
	\item[3.]  Boundary conditions:
	\begin{eqnarray*} &&\hspace{-1cm} (p(0), -p(1)) \in N_{C_0}(x^*(0))\hsm\times\hsm \{0\}+ (0,\lambda_0\partial g(x^*(1))).
	\end{eqnarray*}
	\item[4.] Maximum condition: $ u^*(t)$ maximizes on $\Omega(x(t))$, $[0,1]$-a.e., the map $ u\to  \left\langle p(t),f(x^*(t),u) \right\rangle.$
	\item[5.] $  \nu(t)\nabla_u h(x^*(t),u^*(t))  = p(t)D_uf (x^*(t), u^*(t)) $ for a.a. $t\in [0,1]$.
	\item[6.] $ |\nu(t)|\leq \kappa(\lambda_0+|p(t)|) \text{ for  a.a. } t\in [0,1].$
\end{itemize}
\end{theorem}
We shall postpone the proof of Theorem \ref{theorem PMP regular} to section \ref{section:proof}.

\section{Non-regular mixed constraints} \label{section:PMP non-regular mixed constraint}

In this section, we reformulate problem $(P)$ (see section \ref{sec:introduction}), by rewriting the sweeping term as a pure state constraint (\ref{state constraint refomulated pb}) and a non-regular mixed equality constraint (\ref{mixed non regular constraint refomulated pb}). This results in a non-regular pair, regardless of the regularity of (\ref{mixed constraint}), since at least one of them is non-regular. Due to this non-regularity, purely finitely additive set functions appear as multipliers for the mixed constraints, instead of measures. Moreover, this causes us to lose information provided by the maximization condition.

We begin by giving some definitions related to the notion of purely finitely additive set functions, or the so-called {\it charges}, taken from \cite{rao1983theory}.

On the interval $[0, 1]$, we denote by $\Sigma$ the $\sigma$-algebra of Lebesgue measurable subsets of $[0, 1]$, and by $\mathcal B$ the Borel $\sigma$-algebra in $[0, 1]$.

\begin{definition}
A real-valued finitely additive function $\zeta(.)$ on $\Sigma$ such that $\zeta(\emptyset)=0$ is called a {\it charge}. A charge $\zeta$ is said to be {\it bounded} if $\sup\left\{|\zeta(E)|: E\in\Sigma\right\}<\infty$. The set of all bounded charges is denoted by $ba(\left[0,1\right],\Sigma)$.
\begin{itemize}
    

\item The total variation of $\zeta$, denoted by $|\zeta|$, is the charge defined by
$$|\zeta|(E) = \sup\sum_{i=1}^n
|\zeta (F_i)|$$where the supremum is taken over all finite partitions $\left\{F_i\right\}_{i=1}^n\subset\Sigma$ of $E$.

\item A positive charge $0\leq\zeta$ is called {\it pure} if the only measure $\mu$ such that $0\leq\mu\leq\zeta$ is the trivial measure $\mu=0$. An arbitrary charge $\zeta$ is called pure if $|\zeta|$
is pure.

\item For a Lebesgue measurable function $h : \left[0, 1\right]\to\R$ we write $h=0\ \zeta$-a.e. if,
for all $\epsilon> 0$,
$$|\zeta|(\left\{t\in\left[0, 1\right] : |h(t)| > \epsilon\right\}) = 0.$$
\item The bounded charge $\lambda$ is {\it weakly absolutely continuous} with respect to $\zeta$, and write $\lambda \ll_{\omega} \zeta$, if $E \in \Sigma$, and $|\zeta|(E)=0$ implies that $\lambda(E)=0$. We denote therefore by $ba(\left[0,1\right],\Sigma,\zeta)$ the set of {\it bounded weakly absolutely continuous charges} with respect to $\zeta$.

\item For a function $h\in L_\infty\left(\left[0, 1\right];\R\right)$ and a charge $\zeta\in ba(\left[0,1\right],\Sigma,\ell)$, ($\ell$ stands for the Lebesgue measure) the charge $\Theta$ defined by $d\Theta=h(t)d\zeta$ is given by $\Theta(E):=\int_Eh(t)d\zeta,\; \hbox{for all }E\in\Sigma.$
\end{itemize}
\end{definition}

In the remainder of this section, and in order to adopt the approach used in \cite{Becerriletal2021}, we shall consider the following:
\begin{itemize}
    \item an objective function in the Lagrangian form, instead of the cost $g(.)$ (taken in problem $(P)$);
    \item a truncated (bounded) normal cone $N_C(x)\cap \textrm{int} \B(0,\rho)$, instead of the full one $N_C(x)$. Here $\rho\geq M$ and $\B(0,\rho)$ is the closed ball centered at $0$ of radius $\rho$. 
\end{itemize}
We are interested in the following ``non-regular'' problem:
\begin{align}
		(NP)\mbox{ Minimize} &  \int_0^1L(t,x(t), u(t))dt\nonumber\\ \nonumber
		\mbox{s.t.}\quad &    \dot x(t)= f(x(t),u(t)) - v(t)\nabla_x \psi(x(t))  \;\textrm{a.e.}\nonumber \\
		&  \psi(x(t)) \le 0 \;\forall t\in [0,1] \nonumber\\
	& v(t)\psi(x(t)) =0\; \textrm{for a.a. } t\in[0,1],\nonumber \\
	& (v,u)\in {\cal V}\times \mathcal{U}\nonumber \\
	&  x(0) \in C_0 \subset C\;\; \nonumber \\
	&  h(x(t),u(t))\leq 0\; \textrm{for a.a. } t\in[0,1]\nonumber
	\end{align}
	 where ${\cal V}=\{v\in L_{\infty}([0,1];\R): v(t)\ge 0\;  \forall t\in \left[0,1\right]\}$ and $\mathcal U := L_{\infty}([0,1];\R^m)$, and for the truncated normal cone, we have $|v(t)\nabla_x \psi(x(t))| \leq \rho.$ Here, we take $L:\left[0,1\right]\times\R^n\times\R^m\to\R$ to be a $C^1$ function, and the mixed constraint (\ref{mixed constraint}) to be either regular or non-regular.
	 
	 \begin{definition}\label{def: weak minimize}
	 The feasible process $(x^*,u^*,v^*)$ is called a weak local minimizer to $(NP)$ if there exists $\epsilon >0$ such that $$\int_0^1L(t,x^*(t), u^*(t))dt \le \int_0^1L(t,x(t), u(t))dt$$ holds for all admissible
processes $(x,u,v)$ satisfying
\[ \| x^*-x\|_C < \epsilon, \;\; \textrm{ and }\| (u^*,v^*)-(u,v)\|_{L_\infty} < \epsilon. \]
	 \end{definition}
	 
	 Observe that if $(x^\ast,u^\ast, v^\ast)$ is a solution to $(NP)$, then $(x^\ast, u^\ast)$ must be a solution to $(P)$.
To simplify notation in this section, for a function $q$ depending on $(t,x,u,v)$ we denote by  $q_x$, $q_u$ and $q_v$ its derivatives w.r.t. $x$, $u$ and $v$ respectively, and by $q^\ast(t)$ the function defined by $q(t,x^\ast(t),u^\ast(t), v^\ast(t))$.

Now, we define 
\begin{align*}
  \beta_1(x,u,v):=-v, \beta_2(x,u,v):=\psi(x),
  \beta_3(x,u,v):=h(x,u). 
\end{align*}

Consider the functional:
\begin{align*}
{\cal J}(x,u,v):=&\int_0^1 L(t,x,u)dt,
\end{align*}
and the operators:
\begin{flalign*}
{\cal K}(x,u,v)(t):= & x(t)-x_0-\int_0^t\left[f(x(s),u(s))\right]ds&&\\
&+\int_0^t\left[v(s)\nabla_x\psi(x(s))\right]ds,&&\\
{\cal H}(x,u,v)(t):= & v(t)\psi(x(t)),&&
\end{flalign*}
and ${\cal B}(x,u,v)(t)$ such that, for $j=1,2,3$, ${\cal B}_j(x,u,v)(t):=\beta_j(x(t),u(t),v(t)).$
By our assumptions, we have that these functionals are strictly Fr\'echet differentiable and their derivatives at $\mathbf{x^\ast}:=(x^\ast,u^\ast,v^\ast)$ in the direction of $\mathbf{y}:=(y,b, \kappa)$ are given by
\begin{flalign*}
{\cal J'}(\mathbf{x^\ast};\mathbf{y})&=\int_0^1 \left\{L_x^\ast(s)y(s)+L_u^\ast(s)b(s)\right\}ds,&&\\
{\cal K'}(\mathbf{x^\ast};\mathbf{y})(t)&=y(t)-\int_0^t\left\{f_x^\ast(s)y(s)+f_u^\ast(s)b(s)\right\}ds&&\\
&\hspace{-.8cm}+\int_0^t \left\{v^\ast(s)\nabla_x^2\psi(x^\ast(s))y(s)+\kappa(s)\nabla_x\psi(x(s))\right\}ds &&\\
{\cal H'}(\mathbf{x^\ast};\mathbf{y})(t)&= v^\ast(t)\nabla_x\psi(x^\ast(t))y(t)+\kappa(t)\psi(x^\ast(t))&& \\
{\cal B}_j'(\mathbf{x^\ast};\mathbf{y})(t)&=\beta_{jx}^\ast(t)y(t)+\beta_{ju}^\ast(t)b(t)+\beta_{jv}^\ast(t)\kappa(t) .&&
\end{flalign*}
 Now applying Theorem 4.1 of \cite{Becerriletal2021}, we can deduce the conditions in Theorem \ref{theorem PMP non regular}.


\begin{theorem}\label{theorem PMP non regular}
Suppose that $(x^\ast,u^\ast, v^\ast)$ is a weak local minimizer for $(NP)$. Assume also that the image of $({\cal K'}, {\cal H'})(x^\ast,u^\ast,v^\ast)$ is closed in $C\left(\left[0,1\right],\R^n\right)\times L_\infty\left(\left[0,1\right],\R\right)$. Then, there exist a scalar $\lambda_0\geq 0$, functions $z_j,\  w\in L_1([0, 1]; \R)$ with $z_j\geq 0$, pure charges $\zeta_j$, $\varpi\in ba(\left[0,1\right],\Sigma,\ell)$ with $\zeta_j\geq 0$, ($j=1,2,3$), and functions $\lambda,\ \alpha\in NBV(\left[0,1\right],\R^n)$ satisfying:

\begin{itemize}
\item[a)] Nontriviality condition:
    \begin{equation}\lambda_0+\sum_{j=1}^{3}\left[\left\|z_j\right\|_{L^1}+\zeta_j(\left[0,1\right])\right]+\left\|w\right\|_{L^1}+|\varpi|(\left[0,1\right])>0.\label{nt}\end{equation}
    
    \item[b)] Complementarity slackness condition:
    \[z_j(t)\beta_j^\ast(t)=0\ \hbox{a. e. and }\beta_j^\ast(t)=0\ \zeta_j-\hbox{a.e.}\]
    
    \item[c)] Transversality condition $\lambda(1)=0$ and the costate equation: 
           \begin{flalign*}&-dp(t)^T=(p(t)\hsm-\hsm\alpha(t))^T(f_x^\ast(t)\hsm-\hsm v^\ast(t)\nabla_x^2\psi(x^\ast(t)))dt&&\nonumber\\    &+\lambda_0L_x^\ast(t)dt +(z_2+wv^\ast)(t)\nabla_x\psi(x^\ast(t))dt\\ & \nonumber+z_3(t)h_x^\ast(t)dt&&\end{flalign*}
    where, $p:=\lambda+\alpha$ is absolutely continuous  and, for $(l=1,\ldots, n)$
    \[\alpha_l(0)=-\Theta_l\left[0,1\right],\ \hbox{ and}\ \alpha_l(t)=-\lim_{n}\Theta_l\left[t+\frac{1}{n},1\right]\]
    and $\Theta_l$ is the charge defined by 
    \[d\Theta_l=\psi_{x_l}(x^\ast(t))d\zeta_2+h_{x_l}^\ast(t)d\zeta_3+v^\ast(t)\psi_{x_l}(x^\ast(t))d\varpi.\]

    \item[d)] Stationarity conditions with respect to the control pair $(u,v)$: 
    \begin{flalign}
    0=&\sum_{l=1}^{n}\lambda_lf_{lu}^\ast(t)+\lambda_0L_u^\ast(t)+z_3(t)h_u^\ast(t)
    &&\label{s1}\\
    0=&-\sum_{l=1}^{n}\lambda_l \psi_{x_l}(x^\ast(t))-z_1(t)+w(t)\psi(x^\ast(t))&&\label{s2}\\
    0=&h_u^\ast(t)d\zeta_3&&\label{s3}\\
    0=&-d\zeta_1+\psi(x^\ast(t))d\varpi.\label{s4}&&
    \end{flalign}
\end{itemize}
\end{theorem}

\begin{remark}
\begin{itemize}
    \item Observe that
    from \eqref{s2} and \eqref{s4} we have
$z_1(t)= - \lambda^T(t)\nabla\psi(x^\ast(t))+w(t)\psi(x^\ast(t))$ and 
$d\zeta_1=\psi(x^\ast(t))d\varpi$,
then, we can deduce that the multipliers $z_1$ and $\zeta_1$  can be removed from the nontriviality condition \eqref{nt}. 
    \item We notice that, with this reformulation, the necessary optimality conditions will be satisfied for every admissible process $(\tilde{x},\tilde{u},\tilde{v})$ satisfying $\tilde{\beta_i}(0)=0$ and $\tilde{\beta}_{iu}(0)=0$, for some $i=1,2,3$, considering the charge corresponding to this constraint as the one concentrated in $t=0$ and all the other multipliers as zero.   
    \item The statement of Theorem \ref{theorem PMP non regular} can be extended to the case when the control constraint set is not the whole set $\R^m$. Indeed, expressing the control set $U$ as
\begin{equation*}U:=\left\{u\in\R^m\mid \varphi(u)\leq 0\right\} \label{def: control set non regular}\end{equation*}
where $\varphi:\R^m\to\R^r$, is $C^1$ and the admissible controls $u$ are in $L_\infty([0,1], U)$, then by considering the extra terms
$\beta_{i+3}(x,u,v):=\varphi_i(u), \text{ for } i=1,\ldots,r$,
the stationarity conditions (\ref{s1}), (\ref{s3}) with respect to the control $u$ will take the following form:
 \begin{flalign*}
    0=&\sum_{l=1}^{n}\lambda_lf_{lu}^\ast(t)+\lambda_0L_u^\ast(t)+z_3(t)h_u^\ast(t)&&\nonumber\\
    &+\sum_{i=1}^{r}z_{i+3}(t)\varphi_{iu}^\ast(t)&&\\
    0=&\sum_{i=1}^{r}\varphi_{iu}^\ast(t)d\zeta_{i+3}+h_u^\ast(t)d\zeta_3.&&\label{s3}
    \end{flalign*}
    
\end{itemize}

\end{remark}

\section{Proof of Theorem \ref{theorem PMP regular}} \label{section:proof}
 The proof can be divided into four steps, inspired by \cite{khalil2022maximum}, but adapted to handle the mixed constraint, and utilizing the sweeping process approximation from \cite{de2019optimal} and \cite{zeidan2019sweeping}, as well as the mixed constrained maximum principle from \cite[Theorem 3.1]{arutyunov2016investigation}.
\begin{itemize}
\item[1)] {\bf Sweeping process approximation and solution existence.}

Denote by $(x,u)$ any feasible control process for $(P)$, such that $x(0)=x_0$. Following the approximation idea in \cite{de2019optimal}, we consider an exponential function (Lipschitz in $x$) approximating the non-Lipschitz dynamics (i.e. the sweeping term) of $(P)$.

The approximating system, denoted by $(S_k)$, is:
$$\begin{cases}
\dot x (t)= f(x(t),u(t)) - \gamma_k e^{\gamma_k \psi(x(t))} \nabla_x \psi(x(t)) \;\; \text{a.e.} \\ x(0) = x_{0,k}, \qquad h (x(t),u(t))\leq \delta_k \quad \forall\, t\in [0,1],
\end{cases}$$ where $\{\gamma_k\}$ is a sequence s.t. $\gamma_k \to +\infty$, and $\forall\, k$, $\gamma_k \hsmm\ge\hsmm \frac{2M}{\eta}$ ($M$ and $\eta$ as defined in assumptions H1 and H3), and let $x_{0,k}\hsm\in\hsm C_0$ s.t. 
$x_{0,k}\to x_0$. Moreover, we take $\delta_k =\max_{t\in[0,1]}\hsmm\{h(\bar x_k(t),u(t))\}$ with $\bar x_k$ being the solution to the dynamics in $(S_k)$ for the pair $(x_0, u)$. Clearly, taking $(x(t),u(t))$ such that $\{t\in [0,1]: h(x(t),u(t))=0\}$ is a non-zero measure set, then at the limit, we can guarantee that $\delta_k\to 0$.

Under our assumptions, by using a variation of the arguments in \cite[Lemma 1]{de2019optimal}, we assert the existence of a sequence $\{x_k\}$ also satisfying the mixed constraint (\ref{mixed constraint}). Then, a subsequence converging uniformly to a function $ x\hsm\in\hsm AC([0,1];\R^n)$ can be extracted, and this limit is the unique solution to $(S)\begin{cases}
\dot x (t) \in  f(x(t),u(t)) -N_C(x(t)) \;\; \textrm{a.e.} \\  x(0)=x_0,\qquad h (x(t),u(t))\leq 0\quad [0,1]-a.e.
\end{cases}$

\noindent
Now consider the approximating problem:\vsm
\begin{equation*} ( P_k)\mbox{ Minimize }  g(x(1))\mbox{ subject to } (\bar S_k) \end{equation*} \mbox{ where } $(\bar S_k)$ \mbox{ is the system } 
\begin{eqnarray*}\hspace{0cm} &&  \left\{\begin{array}{l}\hspace{-.1cm}\dot x(t)\hsm= \hsm  f(x(t),u(t))\hsmm -\hsmm\gamma_k e^{\gamma_k \psi(x(t))} \nabla \psi(x(t)) \textrm{ a.e. }\vspace{.1cm}\\
	 x(0) \in C_0 \textrm{ and }   h(x(t),u(t))\leq \delta_k\; \; \forall t\in[0,1] \nonumber\end{array}\right.
\end{eqnarray*}
where $\delta_k $ defined as before but now for data $(x^*(0), u^*)$.

By applying a modified version of the arguments presented in \cite{zeidan2019sweeping}, we can construct a sequence $(\bar x_k,\bar u_k)$ of optimal solutions for $(P_k)$ approximating the solution to $(P)$. 
\vskip2ex
\item[2)] {\bf Ekeland’s Variational Principle.}

By using the arguments above, we can extract a subsequence such that $\bar x_k\to\bar x$ uniformly, where $(\bar x,\bar u)$  is feasible for $(S)$. Indeed, $g(x^*(1))\leq g(\bar x(1))$. We denote by $\tilde x_k$ the unique trajectory for $(\bar S_k)$ with data $(x^*(0),u^*)$. We can directly deduce that $g(\bar x_k(1))\leq g(\tilde x_k(1))$, and $\lim_{k\to \infty}\tilde x_k(1)\hsmm =\hsmm x^*(1)$. These relations guarantee that $ \exists \,\{\varepsilon_k\} \downarrow 0$, such that, for a subsequence (we do not relabel), $g(\tilde x_k(1))\hsmm \leq \hsmm g(\bar x_k(1))\hsmm +\hsmm \varepsilon_k$.

We denote now by $\Delta$ the Ekeland's distance function defined by $$ \Delta( (a,u),(b, v)):=|a-b| +\|u-v\|_{L_1}.$$

For any $\alpha > 0$, $(\tilde x_k,u^*)$ is an $\varepsilon_k$-minimizer to $( P^\alpha_k)$ obtained from $(P_k)$ by
replacing its cost functional by $$g(x(1))\hsmm +\hsmm\alpha\Delta((x(0),u),(x^*(0),u^*)).$$ 
Applying the Ekeland's variational principle, \cite{ekeland1974variational}, there exists a pair $(\hat x_k(0), \hat u_k)$ satisfying $\Delta( (\hat x_k(0),\hat u_k),(x^*_0, u^*)) \leq \sqrt{\varepsilon_k},$
while being the unique solution to $(\bar P^\alpha_k)$ expressed below:
\begin{eqnarray}
\hsmm\mbox{ Min } &&g(x(1))\hsmm +\hsmm(\alpha\hsm +\hsm\sqrt{\varepsilon_k})\Delta((x(0),u),(\hat x_k(0), \hat u_k)) \nonumber\\
&&\hspace{-.7cm}\forall \, (x(0),u)\mbox{ feasible for } (\bar S_k). \vsm\nonumber
\end{eqnarray}
As a consequence, the sequence $\{(\hat x_k(0), \hat u_k)\}$ satisfies $ \| \hat u_k- u^*\|_{L_1}\hsmm\to\hsmm 0 $, and $ \hat x_k\hsmm \to\hsmm x^* $ uniformly.

 \vskip2ex   
    \item[3)] {\bf Maximum Principle to the perturbed problem.}
    
    Under the assumptions on $(P)$ and inherited by $(\bar P^\alpha _k)$, we apply the Maximum Principle in \cite[Theorem 3.1]{arutyunov2016investigation}, to the control process $(\hat x_k, \hat u_k) $, solution to $(\bar P^\alpha_k)$. We define its Pontryagin-Hamilton function as follows:
\begin{eqnarray*}
H_k(x,u,p,\lambda)\hsmm && := \bigl\langle p, f(x,u)\nonumber -\gamma_ke^{\gamma_k\psi(x)} \hsmm\nabla_x\psi(x)\bigr\rangle\hsm  \\
&& -  \lambda\alpha_k|u- \hat u_k(t)|, \textrm{ with } \alpha_k\hsmm =\hsmm \alpha + \sqrt{\varepsilon_k}.
\end{eqnarray*}
Let $(\hat x_k,\hat u_k)$ be an optimal solution to $(\bar P^\alpha_k)$. Then, there exists a set of  multipliers $(p_k,\nu_k,\lambda_k)$ (their dependence on $\alpha_k$ is omited for simplicity) with $p_k \hsm\in\hsm AC([0,1];\R^n)$, a nonnegative 
function $\nu_k\hsm\in\hsm L_\infty([0,1];\R)$, a constant $\kappa >0$, and a number $\lambda_k\geq 0$, satisfying the following conditions:\vskip.2ex
\begin{itemize}
\item[1.] Nontriviality: $ \|p_k\|_{L_\infty}+\lambda_k\neq 0.$
\item[2.] Adjoint equation:
\begin{eqnarray}\nonumber
	-\dot p_k(t) & = \nabla_x H_k(\hat x_k(t),\hat u_k(t), p_k(t),\lambda_k) \\ &  -\nu_k(t)\nabla_x h(\hat x_k(t),\hat u_k(t)).\label{adjoint}\end{eqnarray}
\item[3.] Boundary conditions:
	\begin{eqnarray*} && p_k(0) \in  \lambda_k \alpha_k \xi_{0,k} + N_{C_0}(\hat x_k(0)) \\ &&  -p_k(1) \in \lambda_k\partial g(\hat x_k(1)).\end{eqnarray*}
for some $\xi_{0,k}\hsm \in\hsmm \partial_x |x-\hat x_k(0)|_{|x=\hat x_k(0)}$.
\item[4.] Maximum condition:

$\hat u_k(t)$ maximizes on $\Omega(x)$, $[0,1]$-a.e., the mapping $u\to  H_k(\hat x_k(t),u, p_k(t),\lambda_k) .$

\item[5.] For a.a. $t\in[0,1]$	\begin{eqnarray*}   \nu_k(t)\nabla_u h(x_k(t), u_k(t)) = p_k(t)D_uf (x_k(t), u_k(t))  .	\end{eqnarray*}
	\item[6.] $ |\nu_k(t)|\leq \kappa(\lambda_k+|p_k(t)|) \text{ for  a.a. } t\in [0,1].$
\end{itemize}
\vskip1ex
Observe that omitting the time dependence, and denoting by $f_k^*$, and $f_k^*(u)$ the quantities $f(\hat x_k,\hat u_k)$ and $f(\hat x_k,u)$, respectively. The adjoint equation (\ref{adjoint}) is written as:
\begin{eqnarray*}
-\dot p_k=p_k\Bigl(\hsm D_x f_k^*\hsm-\hsm\gamma_k e^{\gamma_k\psi_k^*} \bigl(\nabla^2_x\psi_k^*+\gamma_kQ_k^*\bigr)\hsm\Bigr)\hsm-\hsm\nu_k\nabla_xh_k^*
\end{eqnarray*}
    where $Q_k^*:= \nabla_x\psi^*_k\otimes\nabla_x\psi^*_k.$

    \vskip2ex
    \item[4)] {\bf Passing to the limit.} 
    
    We will use arguments similar to those in \cite{khalil2022maximum}, with slight modifications to account for the presence of mixed constraints instead of state constraints. We summarize these arguments below without going into technical details.
    
    By the boundedness of the sequence $\{\lambda_k\}$, we can extract a subsequence such that $\lambda_k \to\lambda_0$ for some $\lambda_0 \ge 0 $. Now, given the uniform convergence of $\{\hat x_k\}$, by Gronwall's Lemma, $\{p_k\}$ is uniformly bounded. Using a variation of the arguments in \cite[Theorem 2]{de2019optimal} to accommodate the mixed constraint term, we can prove that $\{\dot{p}_k\}$ is uniformly bounded in $L_1$. Then, there exists a subsequence converging pointwise to a function of bounded variation $p$ (using Helly's theorem). From Helly's second theorem,  $\forall\, \phi\in C([0,1];\R^n)$, $\lim_{k\to\infty}\int_{[0,1]}\phi(t)\dot{p}_k(t)dt= \int_{[0,1]}\phi(t)dp(t)$ for some $p\in BV([0,1];\R^n)$. Moreover, $\|\nu_k\|_{L_\infty}$ is bounded by condition 6., and by a subsequence extraction (not relabeled), there exists $\nu$ such that $\nu_k\to \nu$ in the weak* topology. Now, by taking into consideration equation \eqref{adjoint}, we have the weak convergence in $L_1$ of a subsequence of:
\begin{itemize}
\item[(i)] $\{p_k D_x f_k^*\}$ to $pD_x f^*$;
\item[(ii)] $ \{ \gamma_k e^{\gamma_k\psi_k^*}p_k\nabla^2\psi_k^*\}$ to $ \xi p\nabla^2 \psi^* $. This is a consequence of the uniform bound in $L_2$ of $\{\gamma_k e^{\gamma_k\psi_k^*}\}$ which converges weakly to some $\xi\in L_2$ and to the a.e. convergence of $\{p_k \nabla^2\psi_k^*\}$ to $p\nabla^2 \psi^*$. Also, it is straightforward to see that $\xi$ satisfies the properties in Theorem \ref{theorem PMP regular};
\item[(iii)] $\{\nu_k \nabla_x h^*_k\}$ to $\nu \nabla_xh^*$.
\end{itemize}
Now, let $\bar \eta_k$ be such that, $\forall\,\bar\phi\in C([0,1];\R)$, $\langle \bar \eta_k, \bar \phi\rangle\hsm =\hsm\int_0^1\hsm\gamma_k^2 e^{\gamma_k \psi^*_k}\langle p_k,\nabla \psi_k^*\rangle\bar \phi dt$, which can be shown to be uniformly bounded, and, thus, regarded as a signed Radon measure. Therefore, we can extract a subsequence converging weakly$^*$ to a Radon measure $d\bar \eta$ supported on the set $ I_\psi^*$. Since $\{\nabla_x\psi_k^*\}$ is uniformly bounded, we have that $ \nabla_x\psi_k^*d\bar\eta_k \to \nabla_x\psi^* d\bar\eta \quad \text{weakly}^*.$
Now, we can write the adjoint equation as a measure-driven equation
\begin{equation} -d p = p(D_xf^*\hsmm -\hsmm \xi\nabla^2\psi^*)dt\hsmm  -\hsmm pQ^*d\bar\eta -\nu \nabla_xh^*dt. \label{adjoint2-limit} \vsm\end{equation}
Concerning the boundary conditions, by the closure of the graph of the limiting normal cone and the limiting subdifferential, we deduce the following:
\begin{eqnarray} && \hspace{-1cm}(p(0),  - p(1))\in \alpha\xi_0\hsmm +\hsmm N^L_{C_0}(x^*(0))\hsmm\times\hsmm\lambda_0\partial g(x^*(1)) \label{boundary-limit}\end{eqnarray}
for some $\xi_0\hsm\in\hsm\R^n$ with $|\xi_0|\hsmm =\hsmm 1$.

For the maximization condition, we obtain: $u^*(t)$ maximizes $[0,1]$-a.e. in $\Omega(x)$ the mapping $$u \to \langle p(t), f(t,x^*(t), u)\rangle\hsmm -\hsmm \alpha|u- u^*(t)|.$$
To end the proof, we shall recall that the above limiting multipliers $(p,\nu,\eta,\lambda_0)$ depend on the parameter $\alpha$. Arguing as in \cite{zeidan2019sweeping}, we extract their
limits as $\alpha\to0$, and we obtain the optimality conditions stated in Theorem \ref{theorem PMP regular}.

\end{itemize}


\section{Conclusion}
Problems with dynamics in the form of a sweeping process coupled with mixed constraints pose challenges, particularly in the type of multipliers required for the necessary optimality conditions. Non-regularity automatically arises when reformulating the sweeping process as state and mixed constraints, which necessitates the use of charges in the optimality conditions and precludes any information on the maximization condition. Future work will compare the two approaches used in this paper, study sufficient conditions which guarantee the closure of the image of $(\mathcal{K}',\mathcal{H}')$ required for non-regular cases, explore problems with more complex sweeping sets, and design constraint qualifications to prevent the emergence of pure charges and the degeneracy phenomenon.

\bibliography{ifacconf}                           
\end{document}